\documentclass[oneside,english]{amsart}
\usepackage[T1]{fontenc}
\usepackage{textcomp}
\usepackage[latin9]{inputenc}
\usepackage{amstext}
\usepackage{amsthm}
\usepackage{amssymb}

\makeatletter
\numberwithin{equation}{section}
\numberwithin{figure}{section}
\theoremstyle{plain}
\newtheorem{thm}{\protect\theoremname}[section]
\theoremstyle{definition}
\newtheorem{defn}[thm]{\protect\definitionname}
\theoremstyle{plain}
\newtheorem{lem}[thm]{\protect\lemmaname}

\makeatother

\usepackage{babel}
\providecommand{\definitionname}{Definition}
\providecommand{\lemmaname}{Lemma}
\providecommand{\theoremname}{Theorem}

\begin{document}
\title[{\tiny some generalized central sets theorem along phulara's way}]{A Study Of Some Generalized Central Sets Theorem Near Zero Along Phulara's
Way }
\author{Jyotirmoy Poddar, Sujan Pal}
\address{Department of Mathematics, Techno India University, Saltlake Sector
V, Kolkata-700091, West Bengal, India}
\email{\emph{poddarjyotirmoy5@gmail.com}}
\address{Department of Mathematics, University of Kalyani, Kalyani, Nadia-741235,
West Bengal, India}
\email{\emph{sujan2016pal@gmail.com}}
\begin{abstract}
The Central Sets Theorem near zero was originally proved by Hindman
and Leader. Later a version of Central Sets Theorem was proved by
De, Hindman and Strauss known to be the stronger Central Sets Theorem.
Subsequently many other versions of Central Sets Theorem came, among
which Dev Phulara proved the theorem for a sequence of central sets
instead of taking one set. In this paper, we provide a more general
version of the theorem along Dev Phulara's way near zero.
\end{abstract}

\maketitle

\section{Introduction}

Ramsey theory is a very enriched branch of combinatorics which deals
with the question that when a set with some particular structure,
is partitioned into finitely many cells, then whether one of the cells
also have that structure. This can be approached through many ways,
like Ergodic theory and topological dynamical system, algebra or often
using elementary combinatorics. An age old theorem in this area is
the following, known as the celebrated van der Wearden's theorem,
which we state for a simple motivation.
\begin{thm}[van der Wearden's Theorem\cite{key-10}]
 Let $l,\,r\in\mathbb{N}$. If we have an $r$ colouring of $\mathbb{N}$
then if we are given a length $l$, there exists two numbers $a,\,d\in\mathbb{N}$
so that

\[
\left\{ a,\,a+d,\,\ldots\,,\,a+ld\right\} 
\]

is monochromatic.
\end{thm}

After a long period, in late seventies, Furstenberg, in his famous
work introduced a very interesting set dynamically, known as the central
set \cite{key-3} and proved that central sets have rich combinatorial
structure, which is known as the Central Sets Theorem.
\begin{thm}[{Original Central Sets Theorem\cite[Proposition 8.21]{key-3}}]
 Let $C$ be a central subset of $\mathbb{N}$, let $k\in\mathbb{N}$,
and for each $i\in\left\{ 1,2,\ldots,k\right\} $, let $\langle y_{i,n}\rangle_{n=1}^{\infty}$be
a sequence in $\mathbb{Z}$. There exist sequences $\langle a_{n}\rangle_{n=1}^{\infty}$
in $\mathbb{Z}$ and $\langle H_{n}\rangle_{n=1}^{\infty}$ in $\mathcal{P}_{f}\left(\mathbb{N}\right)$
such that

(1) For each $n\in\mathbb{N}$, $\max H_{n}$< $\min H_{n+1}$ , and

(2) For each $i\in\left\{ 1,2,\ldots,k\right\} $, and $F\in\mathcal{P}_{f}\left(\mathbb{N}\right)$,
we have
\[
\sum_{n\in F}\left(a_{n}+\sum_{t\in H_{n}}y_{i,t}\right)\in C.
\]
\end{thm}

Later in nineties, Hindman and Bergelson took a giant leap in this
area by studying central sets with the help of algebraic structure
of Stone-\v{C}ech compactification of natural numbers \cite{key-11},
which we will discuss in the next section. One of the advantages of
the definition by Bergelson and Hindman over the earlier one is that
from this definition it can be stated easily that a super set of a
central set is again central. From the origin of Furstenberg's Central
Sets Theorem following extensions of Central Sets Theorem has been
established. Hindman, Maleki and Strauss proved a version of Central
Sets Theorem in \cite{key-6}, taking countably infinitely many sequences
and De, Hindman and Strauss proved a version of Central Sets Theorem
in \cite{key-1-1}, taking all sequences at a time. The key focus
of our work is also in this direction, where we will provide several
versions of Central Sets Theorem. 

In case of $\beta\mathbb{N}$ it may be observed that idempotents
live only at infinity, but if we turn our attention to dense subsemigroups
of $\left(\mathbb{R},+\right)$ then idempotents live also near $0$.
The idea first appeared in \cite{key-5}. They showed that there are
localized minimal idempotents near $0$ all of whose members satisfy
some localized Central Sets Theorem conclusion.
\begin{defn}
Let $S$ be a dense subsemigroup of $\left(\left(0,\infty\right),+\right)$.
Then we define the following,

$\mathcal{Z}=\left\{ \begin{array}{c}
\langle\langle y_{i,t}\rangle_{t=1}^{\infty}\rangle_{i=1}^{\infty}\mid\text{for\,each}\,i\in\mathbb{N},\langle y_{i,t}\rangle_{t=1}^{\infty}\,\\
\text{is\,a\,sequence\,in}\,S\cup-S\cup\left\{ 0\right\} \,\text{and\,}\sum_{t=1}^{\infty}\mid y_{i,t}\mid\,\text{converges}\,
\end{array}\right\} $.
\end{defn}

\begin{thm}[Central Sets Theorem near 0]
 Let $S$ be a dense subsemigroup of $\left(\left(0,\infty\right),+\right)$
and let $A$ be a central set near $0$. If we take $Y=\langle\langle y_{i,t}\rangle_{t=1}^{\infty}\rangle_{i=1}^{\infty}\in\mathcal{Z}$,
then there exist sequences $\langle a_{n}\rangle_{n=1}^{\infty}$
in $S$ and $\langle H_{n}\rangle_{n=1}^{\infty}$ in $\mathcal{P}_{f}\left(\mathbb{N}\right)$
such that 

(a) for each $n\in\mathbb{N}$, $a_{n}<\frac{1}{n}$ and $max\,H_{n}<min\,H_{n+1}$
and 

(b) for each $f\in\left\{ f\mid f:\mathbb{N}\rightarrow\mathbb{N}\,\text{and\,for\,all}\,n\in\mathbb{N},f\left(n\right)\leq n\right\} $,
\[
FS\left(\langle a_{n}+\sum_{t\in H_{n}}y_{f\left(n\right),t}\rangle_{n=1}^{\infty}\right)\subseteq A.
\]
\end{thm}

\begin{proof}
\cite[Theorem 4.11]{key-5}.
\end{proof}
In \cite{key-13} the authors extended this theorem for general semigroups.
The above idea has been generalized in \cite{key-9} to get notions
of largeness with respect to filters. Following \cite{key-9} a version
of Central Sets Theorem was proved in \cite{key-4}.

In recent days Dev Phulara provided a very much generalized form of
the Central Sets Theorem in \cite{key-8}, not merely for a single
central set but for a sequence of central sets. This work extended
the work of De, Hindman and Strauss of \cite{key-1-1}, and this is
the key motivation of our work and so we mention that theorem here.
Some notations used in the following theorem are new and will be defined
in the next section.
\begin{thm}
Let $\left(S,\cdot\right)$ be a semigroup and $r$ be an idempotent
in $J\left(S\right)$ and let $\langle C_{n}\rangle_{n=1}^{\infty}$
be a sequence of members of $r$. Then there exist functions 

$m:\mathcal{P}_{f}\left(^{\mathbb{N}}S\right)\rightarrow\mathbb{N}$,
$\alpha\in\times_{F\in\mathcal{P}_{f}\left(^{\mathbb{N}}S\right)}S^{m(F)+1}$
and $\tau\in\times_{F\in\mathcal{P}_{f}\left(^{\mathbb{N}}S\right)}\mathcal{J}_{m(F)}$
such that

(1) if $F,G\in\mathcal{P}_{f}\left(^{\mathbb{N}}S\right)$ and $G\subsetneq F$,
then $\tau\left(G\right)\left(m\left(G\right)\right)$$<$ $\tau\left(F\right)\left(1\right)$
and

(2) when $n\in\mathbb{N}$, $G_{1},G_{2},\ldots,G_{n}\in\mathcal{P}_{f}\left(^{\mathbb{N}}S\right)$,
$G_{1}\subsetneq G_{2}\subsetneq\ldots\subsetneq G_{n}$, and for
each $i\in\left\{ 1,2,\ldots,n\right\} $, $f_{i}\in G_{i}$, and
$l=\mid G_{1}\mid$ one has 
\[
\prod_{i=1}^{n}\left(\left(\prod_{j=1}^{m\left(G_{i}\right)}\alpha\left(G_{i}\right)\left(j\right)\cdot f_{i}\left(\tau\left(G_{i}\right)\left(j\right)\right)\right)\cdot\alpha\left(G_{i}\right)\left(m\left(G_{i}\right)+1\right)\right)\in C_{l}.
\]
\end{thm}

\begin{proof}
\cite[Theorem 3.6]{key-8}.
\end{proof}
The importance of Phulara's version over the previous ones is that
for all the members of the sequence $\langle C_{n}\rangle_{n=1}^{\infty}$
of members of the idempotent $r$ there exists single 
\[
m:\mathcal{P}_{f}\left(^{\mathbb{N}}S\right)\rightarrow\mathbb{N},\alpha\in\times_{F\in\mathcal{P}_{f}\left(^{\mathbb{N}}S\right)}S^{m(F)+1},\tau\in\times_{F\in\mathcal{P}_{f}\left(^{\mathbb{N}}S\right)}\mathcal{J}_{m(F)}.
\]
 But we have to pay a price for that, as the sequence $G_{1}\subsetneq G_{2}\subsetneq\ldots\subsetneq G_{n}$
of\textbf{ $\mathcal{P}_{f}\left(^{\mathbb{N}}S\right)$} starts from
the position of the central set of consideration. 

The paper has been organised as follows. In the next section we are
going to discuss many important definitions and required preliminary
ideas for the understanding of the paper, namely some algebra of the
Stone-\v{C}ech compactification and concepts of idempotents near
0. In section 3 we will prove the Central Sets Theorem near zero along
the idea of Phulara and a consequence.

\section{Definitions and Preliminaries}

We now give a brief review about the Stone-\v{C}ech compactification
of a discrete semigroup. Let $\left(S,\cdot\right)$ be any discrete
semigroup and denote its Stone-\v{C}ech compactification by $\beta S$.
$\beta S$ is the set of all ultrafilters on $S$, where the points
of $S$ are identified with the principal ultrafilters. The basis
for the topology is $\left\{ \bar{A}:A\subseteq S\right\} $, where
$\bar{A}=\left\{ p\in\beta S:A\in p\right\} $. The operation of $S$
can be extended to $\beta S$ making $\left(\beta S,\cdot\right)$
a compact, right topological semigroup with $S$ contained in its
topological center. That is, for all $p\in\beta S$, the function
$\rho_{p}:\beta S\rightarrow\beta S$ is continuous, where $\rho_{p}\left(q\right)=q\cdot p$
and for all $x\in S$, the function $\lambda_{x}:\beta S\rightarrow\beta S$
is continuous, where $\lambda_{x}\left(q\right)=x\cdot q$. For $p,q\in\beta S$
and $A\subseteq S$, $A\in p\cdot q$ if and only if $\left\{ x\in S:x^{-1}A\in q\right\} \in p$,
where $x^{-1}A=\left\{ y\in S:x\cdot y\in A\right\} $. 

Since $\beta S$ is a compact Hausdorff right topological semigroup,
it has a smallest two sided ideal denoted by $K\left(\beta S\right)$,
which is the union of all of the minimal right ideals of $S$, as
well as the union of all of the minimal left ideals of $S$. Every
left ideal of $\beta S$ contains a minimal left ideal and every right
ideal of $\beta S$ contains a minimal right ideal. The intersection
of any minimal left ideal and any minimal right ideal is a group,
and any two such groups are isomorphic. Any idempotent $p$ in $\beta S$
is said to be minimal if and only if $p\in K\left(\beta S\right)$.
Though Central sets was defined dynamically, there is an algebraic
counterpart of this definition, was established by V. Bergelson and
N. Hindman in \cite{key-11} as mentioned in the introduction.
\begin{defn}
Let $S$ be a discrete semigroup. Then a subset $A$ of $S$ is called
central if and only if there is some minimal idempotent $p$ such
that $A\in p$.
\end{defn}

In this context we now need to define a few combinatorially rich sets
which arises now and then in Ramsey theory, later we will also give
these definitions in other settings according as our requirement.
\begin{defn}
Let $\left(S,\cdot\right)$ be a semigroup and $A\subseteq S$, then
\end{defn}

\begin{enumerate}
\item The set $A$ is\emph{ thick} if and only if for any $F\in\mathcal{P}_{f}\left(S\right)$,
there exists an element $x\in S$ such that $F\cdot x\subseteq A$.
For example one can see $\bigcup_{n\in\mathbb{N}}\left[2^{n},2^{n}+n\right]$
is a thick set in $\mathbb{N}$.
\item The set $A$ is called syndetic if and only if there exists a finite
subset $G$ of $S$ such that $\bigcup_{t\in G}t^{-1}A=S$. For example
the set of even and odd numbers are both syndetic in $\mathbb{N}$.
\item $\mathcal{T}=S^{\mathbb{N}}=\left\{ f|f:\mathbb{N}\to S\right\} $.
\item For $m\in\mathbb{N}$, $\mathcal{J}_{m}=\left\{ \left(t\left(1\right),\ldots,t\left(m\right)\right)\in\mathbb{N}^{m}:t\left(1\right)<\ldots<t\left(m\right)\right\} .$
\item Given $m\in\mathbb{N}$, $a\in S^{m+1}$, $t\in\mathcal{J}_{m}$ and
$f\in F$, 
\[
x\left(m,a,t,f\right)=\left(\prod_{j=1}^{m}\left(a\left(j\right)\cdot f\left(t\left(j\right)\right)\right)\right)\cdot a\left(m+1\right)
\]
where the terms in the product $\prod$ are arranged in increasing
order.
\item $A\subseteq S$ is called a $J$-set iff for each $F\in\mathcal{P}_{f}\left(\mathcal{T}\right)$,
there exists $m\in\mathbb{N}$, $a\in S^{m+1}$, $t\in\mathcal{J}_{m}$
such that, for each $f\in\mathcal{T}$,
\[
x\left(m,a,t,f\right)\in A.
\]
\item If the semigroup $S$ is commutative, the definition is rather simple.
In that case, a set $A\subseteq S$ is a $J$-set if and only if whenever
$F\in\mathcal{P}_{f}\left(^{\mathbb{N}}S\right)$, there exist $a\in S$
and $H\in\mathcal{P}_{f}\left(\mathbb{N}\right)$, such that for each
$f\in F$, $a+\sum_{t\in H}f(t)\in A$.
\end{enumerate}
It shouldbe noted that a set is thick if it contains a translation
of any finite subset. Also with a finite translation, if the set covers
the entire semigroup, then it will be called a Syndetic set.

We next define the notion of idempotents near zero, originally introduced
by Hindman and Leader in \cite{key-5}.

We work with $S=\left(\left(0,\infty\right),+\right)$. We have been
considering those semigroups which are dense in $S$ with respect
to the natural topology of $S$. When we want to discuss the Stone-\v{C}ech
compactification of such a semigroup $S$, we have to shift to $S_{d}$,
the set $S$ with the discrete topology.
\begin{defn}
Let $S$ be a dense subset of $\left(\left(0,\infty\right),+\right)$.
Then 
\[
0^{+}\left(S\right)=\left\{ p\in\beta S_{d}:\left(\forall\epsilon>0\right)\left(0,\epsilon\right)\cap S\in p\right\} .
\]
\end{defn}

We now have to recall the notions of combinatorially rich sets near
zero from the literature.
\begin{defn}
Let $S$ be a dense subsemigroup of $\left(\left(0,\infty\right),+\right)$.
and let $A\subseteq S$. 
\end{defn}

\begin{enumerate}
\item $A$ is a central set near zero if and only if there exists an idempotent
$p$ in the smallest ideal of $0^{+}\left(S\right)$ with $A\in p$.
\item A subset $A$ of $\left(0,1\right)$ is Syndetic near $0$ if and
only if $\forall\epsilon>0$ there exists $F\in\mathcal{P}_{f}\left(0,\epsilon\right)$
and $\delta>0$ such that $\left(0,\delta\right)\subseteq\bigcup_{t\in F}\left(-t+A\right)$.
\item We say that $f:\mathbb{N}\rightarrow S$ is near zero if $lim\left(f\left(\mathbb{N}\right)\right)=0$.
The collection of all sequences that is near zero is denoted by $\mathcal{T}_{0}$. 
\item A subset $A$ of $\left(0,1\right)$ is called $J$-set near $0$
iff whenever $F\in\mathcal{P}_{f}\left(\mathcal{T}_{0}\right)$ and
$\delta>0$ , there exists $a\in S\cap\left(0,\delta\right)$ and
$H\in\mathcal{P}_{f}\left(\mathbb{N}\right)$ such that for each $f\in F$,
$a+\sum_{t\in H}f\left(t\right)\in A$.
\item $J_{0}\left(S\right)=\left\{ p\in0^{+}:\forall A\in p,A\:\text{is}\:\text{a}\:J\text{-set \,near\,0}\right\} $. 
\end{enumerate}

\section{Central Sets Theorem Near Zero Along Phulara's Way}

In this section we will show that the Central Sets Theorem near $0$
can be modified in the direction of Dev Phulara, i.e we show that
the conclusion of the theorem is true if we take a sequence of central
sets instead of a single central set. But before that we need to state
a lemma first.
\begin{lem}
\label{lemma=000020for=000020J=000020set} Let, $S$ be a dense subsemigroup
of $\left(\left(0,\infty\right),+\right)$ and $A\subseteq S$ is
a $J$\textminus set near zero. Whenever $m\in\mathbb{N}$ and $F\in\mathcal{P}_{f}\left(\mathcal{T}_{0}\right)$
and $\delta>0$, there exist $a\in S\cap\left(0,\delta\right)$ and
$H\in\mathcal{P}_{f}\left(\mathbb{N}\right)$ such that $\text{min}\,H>m$
and for each $f\in F$, $a+\sum_{t\in H}f\left(t\right)\in A$.
\end{lem}

\begin{proof}
\cite[Lemma 3.3]{key-13}.
\end{proof}
Now we are in a position to prove our required version.
\begin{thm}
\label{CST=000020near=000020zero=000020along=000020Phulara} Let $S$
be a dense subsemigroup of $\left(\left(0,\infty\right),+\right)$,
let $p$ be an idempotent in $J_{0}\left(S\right)$, and let $\langle C_{n}\rangle_{n=1}^{\infty}$
be a sequence of members of $p$. Then for each $\delta\in\left(0,1\right),$
there exists $\alpha_{\delta}:\mathcal{P}_{f}\left(\mathcal{T}_{0}\right)\rightarrow S$
and $H_{\delta}:\mathcal{P}_{f}\left(\mathcal{T}_{0}\right)\rightarrow\mathcal{P}_{f}\left(\mathbb{N}\right)$such
that

1. $\alpha_{\delta}\left(F\right)<\delta$ for each $F\in\mathcal{P}_{f}\left(\mathcal{T}_{0}\right)$;

2. if $F,G\in\mathcal{P}_{f}\left(\mathcal{T}_{0}\right),F\subsetneq G$,
then $\text{max}H_{\delta}\left(F\right)<\text{min}H_{\delta}\left(G\right)$
and

3. if $m\in\mathbb{N}$and $G_{1},G_{2},...,G_{m}\in\mathcal{P}_{f}\left(\mathcal{T}_{0}\right),G_{1}\subsetneq G_{2}\subsetneq\ldots\subsetneq G_{m}$,
$f_{i}\in G_{i}$ for each $i=1,2,\ldots,m$, and $|G_{1}|=r$, then
\[
\sum_{i=1}^{m}\left(\alpha_{\delta}\left(G_{i}\right)+\sum_{t\in H_{\delta}\left(G_{i}\right)}f_{i}\left(t\right)\right)\in C_{r}.
\]
\end{thm}

\begin{proof}
We may assume that $\langle C_{n}\rangle_{n=1}^{\infty}$ to be decreasing.
For each $n$, let $C_{n}^{*}=\left\{ x\in C_{n}:-x+C_{n}\in p\right\} $.
Then $C_{n}^{*}\in p$ and by \cite[Corollary 4.14]{key-7}, for each
$x\in C_{n}^{*}$, $-x+C_{n}^{*}\in p$.

Let $\delta\in\left(0,1\right)$ be given. We define $\alpha_{\delta}\in S$
and $H_{\delta}\in\mathcal{P}_{f}\left(\mathbb{N}\right)$ for $F\in\mathcal{P}_{f}\left(\mathcal{T}_{0}\right)$
by induction on $\mid F\mid$ so that

(1) $\alpha_{\delta}\left(F\right)<\delta$.

(2) if $F,G\in\mathcal{P}_{f}\left(\mathcal{T}_{0}\right),G\subsetneq F$,
then $\text{max}H_{\delta}\left(G\right)<\text{min}H_{\delta}\left(F\right).$

(3) if $m\in\mathbb{N}$and $\emptyset\neq G_{1}\subsetneq G_{2}\subsetneq\ldots\subsetneq G_{m}=F$,
$f_{i}\in G_{i},$for each $i=1,2,\ldots,m$, and $|G_{1}|=r$, then
\[
\sum_{i=1}^{m}\left(\alpha_{\delta}\left(G_{i}\right)+\sum_{t\in H_{\delta}\left(G_{i}\right)}f_{i}\left(t\right)\right)\in C_{r}^{\star}.
\]

Let $f\in\mathcal{T}_{0}$ and let $F=\left\{ f\right\} $. Since
$C_{1}\in p$ and $p\in J_{0}\left(S\right)$, $C_{1}^{\star}$ is
J-set near zero, so for given $\delta>0,$pick $a\in S\cap\left(0,\delta\right)$
and $L\in\mathcal{P}_{f}\left(\mathbb{N}\right)$ such that $a+\sum_{t\in L}f\left(t\right)\in C_{1}^{\star}.$

Let $\alpha_{\delta}\left(F\right)=a$ and $H_{\delta}\left(F\right)=L$.
Then thehypotheses are satisfied, (2) is vacuously true.

Now assume that $F\in\mathcal{P}_{f}\left(\mathcal{T}_{0}\right)$,
$\mid F\mid=n>0$, and $\alpha_{\delta}\left(G\right)$ and $H_{\delta}\left(G\right)$
have been defined for all proper subsets $G$ of $F$, satisfying
the induction hypotheses. 

Let $K_{\delta}=\bigcup\left\{ H_{\delta}\left(G\right):\emptyset\neq G\subsetneq F\right\} $
and let $d=\text{max}\,K_{\delta}.$ 

For $r\in\left\{ 1,2,\ldots,n-1\right\} $, let
\[
M_{\delta}^{r}=\left\{ \begin{array}{c}
\sum_{i=1}^{s}\left(\alpha_{\delta}\left(G_{i}\right)+\sum_{t\in H_{\delta}\left(G_{i}\right)}f_{i}\left(t\right)\right):\\
s\in\left\{ 1,2,\ldots,n-1\right\} ,\emptyset\neq G_{1}\subsetneq G_{2}\subsetneq\ldots\subsetneq G_{s}\subsetneq F,\\
f_{i}\in G_{i}\,\text{for}\,i\in\left\{ 1,2,\ldots,s\right\} ,\,\text{and\,}\mid G_{1}\mid=r
\end{array}\right\} .
\]
 Then each$M_{\delta}^{r}$ is finite and by hypothesis (3), $M_{\delta}^{r}\subseteq C_{r}^{\star}$.
Let 
\[
A=C_{n}^{\star}\cap\bigcap_{r=1}^{n-1}\bigcap_{x\in M_{\delta}^{r}}\left(-x+C_{r}^{\star}\right).
\]

Then, $A\in p,$and so $A$ is a $\text{J}$-set near zero. By \ref{lemma=000020for=000020J=000020set}
pick $a\in S\cap\left(0,\delta\right)$ and $L\in\mathcal{P}_{f}\left(\mathbb{N}\right)$
such that $\text{min}\,L>d$ and for each $f\in F$, $a+\sum_{t\in L}f\left(t\right)\in A$.
Let $\alpha_{\delta}\left(F\right)=a$ and $H_{\delta}\left(F\right)=L$.
Hypothesis (1) holds directly and since $\text{min}\,L>d$, hypothesis
(2) holds.

To verify hypothesis (3) let $\emptyset\neq G_{1}\subsetneq G_{2}\subsetneq\ldots\subsetneq G_{m}=F$,
let $f_{1},f_{2},\ldots,f_{m}$ be given such that each $f_{i}\in G_{i}$,
and let $r=\mid G_{1}\mid$. and let for each $i\in\left\{ 1,2,\ldots,s\right\} ,$$f_{i}\in G_{i}$.
Assume first that $m=1$. Then hypothesis (3) holds because $\alpha_{\delta}\left(F\right)+\sum_{t\in H_{\delta}\left(F\right)}f\left(t\right)\in A\subseteq C_{n}^{*}.$

Now assume that $m>1$, and let $r=\mid G_{1}\mid$. Let $x=\sum_{i=1}^{m-1}\left(\alpha_{\delta}\left(G_{i}\right)+\sum_{t\in H_{\delta}\left(G_{i}\right)}f_{i}\left(t\right)\right)$.

Then $x\in M_{\delta}^{r}$, and $\alpha_{\delta}\left(F\right)+\sum_{t\in H_{\delta}\left(F\right)}f_{m}\left(t\right)\in A\subseteq\left(-y+C_{r}^{\star}\right)$
as required.
\end{proof}
We now discuss a combinatorial result which is classic but we generalize
together for near zero and along a sequence of central sets.
\begin{thm}
\label{CST=000020another=000020version} Let $S$ be a dense subsemigroup
of $\left(\left(0,\infty\right),+\right)$. Let $\left(C_{n}\right)_{n=1}^{\infty}$
be a sequence of central subsets belonging to an idempotent $p\in J_{0}\left(S\right)$.
Let $k\in\mathbb{N}$ and for each $l\in\left\{ 1,2,\ldots,k\right\} $
let $\left\langle y_{l,n}\right\rangle _{n=1}^{\infty}$ be a sequence
in $\mathcal{T}_{0}$. Then there exists a sequence $\left\langle a_{n}\right\rangle _{n=1}^{\infty}$
in $\mathcal{T}_{0}$ and a sequence $\left\langle H_{n}\right\rangle _{n=1}^{\infty}$
in $\mathcal{P}_{f}\left(\mathbb{N}\right)$ such that $\text{max}\,H_{n}<\text{min}\,H_{n+1}$
for each $n\in\mathbb{N}$ and such that for each $l\in\left\{ 1,2,\ldots,k\right\} $and
$F\in\mathcal{P}_{f}\left(\mathbb{N}\right)$ with $\text{min}\,F=m$,
we have 
\[
\sum_{n\in F}\left(a_{n}+\sum_{t\in H_{n}}y_{l,t}\right)\in C_{m}.
\]
\end{thm}

\begin{proof}
We may assume that $\left(C_{n}\right)_{n=1}^{\infty}$ is downward
directed. Pick $\alpha_{\delta}$ and $H_{\delta}$ as guaranteed
by \ref{CST=000020near=000020zero=000020along=000020Phulara}. Now
choose $\left\langle \gamma_{u}\right\rangle _{u=1}^{\infty}$in this
manner. Let
\[
\gamma_{1}\in\mathcal{T}_{0}\setminus\left\{ \left\langle y_{1,n}\right\rangle _{n=1}^{\infty},\left\langle y_{2,n}\right\rangle _{n=1}^{\infty},\ldots,\left\langle y_{k,n}\right\rangle _{n=1}^{\infty}\right\} 
\]
 which is defined due to the non-triviality of $S$. For $u\in\mathbb{N},$
pick 
\[
\gamma_{u+1}\in\mathcal{T}_{0}\setminus\left(\left\{ \left\langle y_{1,n}\right\rangle _{n=1}^{\infty},\left\langle y_{2,n}\right\rangle _{n=1}^{\infty},\ldots,\left\langle y_{k,n}\right\rangle _{n=1}^{\infty}\right\} \cup\left\{ \gamma_{1},\gamma_{2},\ldots,\gamma_{u}\right\} \right).
\]
 This choiceis possible due to the fact that $\mathcal{T}_{0}$ is
infinite. Now define

\[
G_{u}=\left\{ \left\langle y_{1,n}\right\rangle _{n=1}^{\infty},\left\langle y_{2,n}\right\rangle _{n=1}^{\infty},\ldots,\left\langle y_{k,n}\right\rangle _{n=1}^{\infty}\right\} \cup\left\{ \gamma_{1},\gamma_{2},\ldots,\gamma_{u}\right\} .
\]
 Let $a_{\delta}^{u}=\alpha_{\delta}\left(G_{u}\right)$ and $H_{\delta}^{u}=H_{\delta}\left(G_{u}\right)$.
Let $l\in\left\{ 1,2,\ldots,k\right\} $ and let $F\in\mathcal{P}_{f}\left(\mathbb{N}\right)$
which is enumerated as $\left\{ n_{1},n_{2},\ldots,n_{s}\right\} $
so that $m=n_{1}$. Then we have 
\[
G_{m}=G_{n_{1}}\subsetneq G_{n_{2}}\subsetneq\ldots\subsetneq G_{n_{s}}.
\]
 Also we have that for each $i\in\left\{ 1,2,\ldots,s\right\} $,
$\left\langle y_{l,t}\right\rangle _{t=1}^{\infty}\in G_{n_{i}},$and
$\mid G_{n_{1}}\mid=m+k$. So we have, 
\[
\sum_{n\in F}\left(a_{\delta}^{n}+\sum_{t\in H_{\delta}^{n}}y_{l,t}\right)=\sum_{i=1}^{s}\left(\alpha_{\delta}\left(G_{n_{i}}\right)+\sum_{t\in H_{\delta}\left(G_{n_{i}}\right)}y_{l,t}\right)\in C_{m+k}\subseteq C_{m}.
\]
\end{proof}
$\vspace{0.3in}$

\textbf{Acknowledgment:} The second author acknowledges the Grant
CSIR-UGC NET fellowship with file No. 09/106(0199)/2019-EMR-I of CSIR-UGC
NET. 

$\vspace{0.3in}$


\begin{thebibliography}{10}
\bibitem[1]{key-13} E. Bayatmanesh, M. Akbari Tootkaboni, Central
sets theorem near zero, Topology and its Applications 210(2016) 70-80.

\bibitem[2]{key-11} V. Bergelson, N. Hindman, \emph{Nonmetrizable
topological dynamics and Ramsey theory}, Trans. Amer. Math. Soc. 320
(1990), 293-320.

\bibitem[3]{key-1-1} D. De, N. Hindman and D. Strauss, \textit{A
new and stronger }\emph{Central Sets Theorem}, Fundamenta Mathematicae
199 (2008), 155-175.

\bibitem[4]{key-2} D. De, N. Hindman and D. Strauss, \emph{Sets central
with respect to certain subsemigroups of $\beta S_{d}$,} Topology
Proceedings 33 (2009), 55-79.

\bibitem[5]{key-3} H. Furstenberg, \textit{Recurrence in ergodic
theory and combinatorial number theory,} Princeton University Press,
Princeton, New Jersey, (1981)

\bibitem[6]{key-4} S. Goswami , J. Poddar, \emph{Central Sets Theorem
along filters and some combinatorial consequences}, Indagationes Mathematicae
33 (2022), 1312-1325.

\bibitem[7]{key-5} N. Hindman, I. Leader{\normalsize , }\emph{The
semigroup of ultrafilters near $0$,}{\normalsize{} }Semigroup Forum
59 (1999), 33-55.

\bibitem[8]{key-6} N. Hindman, A. Maleki and D. Strauss{\normalsize ,
}\emph{Central sets and their combinatorial characterization},{\normalsize{}
}J. Comb. Theory (Series A) 74 (1996), 188-208.

\bibitem[9]{key-7} N.Hindman, D.Strauss{\normalsize , }\emph{Algebra
in the Stone-\v{C}ech Compactifications: theory and applications}{\normalsize ,
}second edition, de Gruyter, Berlin, 2012.

\bibitem[10]{key-8} D.Phulara, \emph{A generalized Central Sets Theorem
and applications}, Topology and its Applications 196(2015) 92-105.

\bibitem[11]{key-9} O. Shuungula,Y. Zelenyuk and Y. Zelenyuk{\normalsize ,
}\emph{The closure of the smallest ideal of an ultrafilter semigroup}{\normalsize ,
}Semigroup Forum 79 (2009) , 531-539.

\bibitem[12]{key-10} B. van der Waerden, \textit{Beweis einer Baudetschen
Vermutung}, Nieuw Arch. Wiskunde 19 (1927), 212-216.

\end{thebibliography}
\end{document}